\documentclass[12pt]{article}

\usepackage{latexsym,amsmath,amsthm,amssymb,fullpage}

\usepackage{color}

\usepackage{latexsym,amsmath,amsthm,amssymb,fullpage}

\usepackage[utf8]{inputenc}

\usepackage{amssymb}
\usepackage{epstopdf}
\usepackage{amsmath}
\usepackage{amsthm}
\usepackage{amsfonts}
\usepackage{epstopdf}

\newtheorem{theorem}{Theorem}
\newtheorem{lemma}[theorem]{Lemma}

\theoremstyle{definition}

\theoremstyle{remark}

\newcommand\beq{\begin {equation}}
\newcommand\eeq{\end {equation}}
\newcommand\beqs{\begin {equation*}}
\newcommand\eeqs{\end {equation*}}

\newcommand\N{{\mathbb N}}

\newcommand\R{{\mathbb R}}

\newcommand\PP{{\mathbb P}}

\newcommand\eu{{\rm e}}

\begin {document}

$$   $$
\begin{center}
\textbf{\LARGE A law of the iterated logarithm\\}
\medskip
\textbf {\LARGE for directed last passage percolation}
\vskip 6mm
\textit{\large Michel Ledoux}\\
\vskip 3mm
{\large University of Toulouse, France}\\
\end{center}

\vskip 10mm

\begin {abstract} 
Let ${\widetilde H}_N$, $N \geq 1$, be the last passage times of directed percolation
on rectangles $[(1,1), ([\gamma N], N])]$ in $\N^2$ over exponential or geometric independent
random variables, rescaled to converge to the Tracy-Widom distribution.
It is proved that for some $\alpha_{\sup} >0$,
$$
\alpha_{\sup} \, \leq \, 
    \limsup_{N \to \infty} \frac {{\widetilde H}_N}{(\log \log N)^{2/3}} \, \leq \, \Big ( \frac 34 \Big)^{2/3}.
$$
with probability one, and that $\alpha_{\sup} = \big ( \frac 34 \big)^{2/3}$ provided
a commonly believed tail bound holds. The result is in contrast with the normalization
$(\log N)^{2/3}$ for the largest eigenvalue of a GUE matrix recently put forward
by E.~Paquette and O.~Zeitouni. The proof relies on sharp tail bounds and superadditivity,
close to the standard law of the iterated logarithm. A weaker result on the liminf is also
discussed.
\end {abstract}

\vskip 4mm

\section {Introduction and main results}

Let ${(X_{i,j})}_{(i, j) \in \N^2}$ be an infinite array of independent exponential random variables
with parameter $1$. For $M \geq N \geq 1$, let
$$
H(M,N) \, = \, \max \bigg \{ \sum_{(i,j) \in \pi} X_{i,j} \, ; \, \,  \pi \in \Pi_{M,N} \bigg \},
$$
where $ \Pi_{M,N}$ is the set of all up/right paths in $\N^2$ joining $(1,1)$ to $(M,N)$,
be the directed last passage time on the rectangle $[(1,1), (M,N)]$ in $\N^2$.

It is a result due to K.~Johansson \cite {J00} that for each $\gamma \geq 1$,
$$
{\widetilde H }_N \, = \, 
\frac { H([\gamma N], N) - a N}{bN^{1/3}} \, ,
$$
where $ a = a(\gamma) = (1 + \sqrt \gamma)^2$
and $b = b (\gamma ) = \gamma^{-1/6} (1 + \sqrt \gamma )^{4/3}$,
converges as $N \to \infty$ to the Tracy-Widom distribution $F_2$. As is by now classical,
the distribution $F_2$ arises as the limit of rescaled largest eigenvalue 
$$
{\widetilde \lambda}_N \, = \,
N^{1/6} \big ( \lambda_N - 2 \sqrt N \, \big)
$$
of the Gaussian Unitary Ensemble (GUE) of size $N$ consisting of an Hermitian
matrix with entries that are independent
(up to the symmetry condition) complex Gaussian variables with mean zero and variance $1$.

In addition to this result, it is also shown in \cite {J00} that $H(M,N)$ has the same distribution
as the largest eigenvalue of the Laguerre Unitary Ensemble, that is of a complex Wishart matrix
$A A^*$ where $A$ is an $N \times M$ matrix with entries that are independent complex Gaussian variables
with mean zero and variance $\frac 12$.

It was recently established by E.~Paquette and O.~Zeitouni \cite {PZ15}
that (whenever the GUE is constructed
from a given infinite array of Gaussian variables on the same probability space),
$$
\limsup_{N \to \infty} \frac {{\widetilde \lambda}_N}{(\log N)^{2/3}} 
    \, = \, \Big ( \frac 14 \Big)^{2/3}
$$
almost surely. It is reasonable to expect that a similar behaviour, of order $(\log N)^{2/3}$,
holds for the largest eigenvalue of a Wishart matrix. One crucial aspect of the investigation
\cite {PZ15} is that the subsequence $N = k^3$ carries much of the almost sure behaviour
in contrast with the standard geometric subsequences in the classical block
argument of the law of the iterated logarithm
(see e.g. \cite {B86} for a survey on the law of the iterated logarithm
and some relevant classical references). The work \cite {PZ15} also presents
a result on the liminf, although with non-optimal limits at this point.

However, in the last passage percolation representation, the almost sure behaviour
actually turns out to be much smaller and of more classical $\log \log$ type.

\begin {theorem} \label {thm.limsup}
There exists $\alpha_{\sup} >0$ such that
$$
\alpha_{\sup} \, \leq \, 
    \limsup_{N \to \infty} \frac {{\widetilde H}_N}{(\log \log N)^{2/3}} \, \leq \, \Big ( \frac 34 \Big)^{2/3}
$$
with probability one.
\end {theorem}

It is expected that $\alpha_{\sup} = \big ( \frac 34 \big)^{2/3}$ and we actually provide a proof of
it based on the suitable tail estimate which is commonly believed to hold true.

There is a similar, although weaker, result for the liminf.

\begin {theorem} \label {thm.liminf}
There exists $0 < \alpha_{\inf} < \infty$ such that
$$
   - \alpha_{\inf}  \, \leq \, 
    \liminf_{N \to \infty} \frac {{\widetilde H}_N}{(\log \log N)^{1/3}} 
$$
with probability one.
\end {theorem}

We have not been able to show the existence of $\beta_{\inf} >0$ such that
$$
     \liminf_{N \to \infty} \frac {{\widetilde H}_N}{(\log \log N)^{1/3}} \, \leq \, - \beta_{\inf}
$$
with probably one. It may be conjectured that $\alpha_{\inf} = \beta _{\inf}= (12)^{1/3}$.

The proofs of Theorem~\ref {thm.limsup} and \ref {thm.liminf} rely on precise
tail inequalities on the distribution of $H([\gamma N], N)$
together with blocking arguments on the path representation. Roughly speaking,
the powers $\frac 23$ and $\frac 13$ reflect the right and left tails of the Tracy-Widom
distribution (cf.~e.g.~\cite {AGZ10})
$$
 1 - F_2(x) \, \sim \, \eu^{- \frac 43  x^{3/2}} \quad \text{as} \quad x \to +\infty, \qquad
F_2(x) \, \sim \, \eu^{-\frac {1}{12}  x^3} \quad \text{as} \quad x \to - \infty, 
$$
whereas the $\log \log$ is the result of a block argument along geometric
subsequences. One main difference with the random matrix models is that
the path representation allows for (point-wise) superadditivity, not available
for extremal eigenvalues, which lead to the almost sure $\log \log$ behaviour.
As a consequence, the proofs here turn out to be
simpler than the study developed in \cite {PZ15} which is making use
of delicate decorrelation estimates obtained via a hard analysis
of the determinantal kernel of the GUE.

The picture on the tail inequalities used in this note is a bit incomplete at this point, 
impacting the main conclusions, although sharp versions should reasonably hold true.

First, the large deviation estimates developed by K.~Johansson in \cite {J00} show that
\beq \label {eq.ldpright}
\lim_{N \to \infty} \frac {1}{N} \log
\PP \big ( H([\gamma N], N) \geq (a + \varepsilon ) N \big)\, = \, - J(\varepsilon) 
\eeq
for each $\varepsilon >0$ where $J$ is an explicit rate function such that $J(x) >0$
if $x>0$. On the left of the mean,
\beq \label {eq.ldpleft}
\lim_{N \to \infty} \frac {1}{N^2} \log
\PP \big ( H([\gamma N], N) \leq (a - \varepsilon ) N \big)\, = \, - I(\varepsilon) 
\eeq
for each $\varepsilon >0$ where $I(x) >0$ if $x>0$.

A superadditivity argument (see \cite {J00} and below) actually
allows in \eqref {eq.ldpright} for the upper bound
\beq \label {eq.rightupperj}
\PP \big ( H([\gamma N], N) \geq (a + \varepsilon ) N \big)\, \leq \, \eu^{ - J(\varepsilon)N}
\eeq
for any $N\geq 1$ and $\varepsilon >0$. The relevant information on $J$ is that (cf.~\cite {J00})
\beq \label {eq.j}
\lim_{\varepsilon \to 0} \frac {J(\varepsilon)}{\varepsilon^{3/2}} \, = \, \frac {4}{3 b^{3/2}} \, .
\eeq
(See also \cite {L07}.)

We will also need a lower bound on the probability in \eqref {eq.rightupperj}, but
the sharp version is not so explicit in the literature (see below).
Nevertheless, in the random matrix interpretation,
we can make use of the results of \cite {LR10} from which
\beq \label {eq.rightlowerlr}
\PP \big ( H([\gamma N], N) \geq (a + \varepsilon ) N \big)\, \geq \, c \, \eu^{ - C \varepsilon ^{3/2}N}
\eeq
for every $0 \leq \varepsilon \leq 1$ and $N \geq 1$, where $c,C >0$ only depend on $\gamma$.
This inequality is actually not detailed in \cite {LR10} but, as explained there, the same
arguments may be used.

Below the mean, following \cite {LR10} in the random matrix description,
for some $c, C >0$ only depending on $\gamma$,
\beq \label {eq.leftupperlr}
\PP \big ( H([\gamma N], N) \leq (a - \varepsilon ) N \big)\, \leq \, C \, \eu^{ - c \varepsilon ^{3}N^2}
\eeq
for every $0 < \varepsilon \leq a$ and $N \geq 1$.

The investigation here may actually be considered similarly for random variables $X_{i,j}$
with a geometric distribution
rather than exponential as in the original contribution \cite {J00},
and Theorems~\ref {thm.limsup} and \ref {thm.liminf} extend to this setting. 
The fluctuations and large deviations are actually established initially for geometric
distributions in \cite {J00} (with suitable values of $a$, $b$ and a suitable $J$ function),
the exponential case being seen as the limit of the geometric model
with parameter tending to $1$. The tail inequality \eqref {eq.rightupperj} holds similarly.
At the level of sharp tail inequalities in the context of geometric random variables,
a refined Riemann-Hilbert analysis on the determinantal structure
of the underlying Meixner Ensemble has been developed in \cite {BDMLMZ01}
to show that below the mean
\beq \label {eq.leftbaik}
\log \PP \big ( H([\gamma N], N) \leq a N  - x b N^{1/3}\big)
   \, = \, - \frac {x^3}{12} + O (x^4 N^{-2/3}) + O(\log x)
\eeq
uniformly over $ M \leq x \leq \delta N^{2/3}$ for some (large) constant $M>0$ and
some (small) constant $\delta >0$, and every $N$ large
enough. Although not written explicitly, it is expected that the same method
(even in a simpler form) may be used above the mean to yield
\beq \label {eq.rightbaik}
\log \PP \big ( H([\gamma N], N) \geq a N  + x b N^{1/3}\big)
   \, = \, - \frac 43 \, x^{3/2} + O (x^2 N^{-1/3}) + O(\log x)
\eeq
uniformly over $ M \leq x \leq \delta  N^{1/3}$ for some (large) constant $M>0$ and
some (small) constant $\delta >0$, and every $N$ large enough.

The same Riemann-Hilbert analysis on the Laguerre Unitary Ensemble
yields \eqref {eq.leftbaik} in the exponential case, and supposedly also
\eqref {eq.rightbaik} (as well as in the GUE setting).
In particular, \eqref {eq.leftbaik} provides a sharp
(two-sided) version of \eqref {eq.leftupperlr} while \eqref {eq.rightbaik} matches
\eqref {eq.rightupperj} and would provide the sharp version
of \eqref {eq.rightlowerlr}. Taking \eqref {eq.rightbaik} for granted, we will prove the sharp version of
Theorem~\ref {thm.limsup} with $\alpha = \big ( \frac 34 \big)^{2/3}$
both in the exponential and geometric cases.

\section {Proofs}

Before addressing the proof of the main results, we emphasize a number of useful tools.
To start with, to avoid some unessential technicalities, in the definition of $H(M,N)$
(and related quantities of the same type), we will actually
consider sums $ \sum_{(i,j) \in \pi} X_{i,j} - X_{1,1}$ (that is omitting the common
initial point of all paths).
It is immediate that this change does not alter any of the limits studied here.

Next, we recall from \cite {J00} the simple but basic superadditivity property.
For simplicity, we write below $W_N = H([\gamma N], N) $, $\gamma \geq 1$
being fixed throughout this work.
Whenever $N \leq L$, let $W_{[N,L]}$ be
the maximum of up/right paths joining $([\gamma N], N)$ to $([\gamma L], L)$ in $\N^2$
(with therefore the preceding convention, that is omitting $X_{[\gamma N], N}$ in the sums). Then,
as is immediate,
\beq \label {eq.superadditivity}
 W_N + W_{[N,L]}  \, \leq  \, W_L.
\eeq

Finally, it will be useful to rely on the following maximal inequality of the type of
the classical Ottaviani inequality for sums of independent random variables
or vectors (cf.~\cite {LT91}).

\begin {lemma} \label {lem.maximal}
For any real numbers $t,s$, and any integers $ 1 \leq K \leq L$,
$$
\PP \Big ( \max_{K \leq N \leq L} ( W_N - a N ) \geq t \Big)
  \, \leq \,  \frac { \PP ( W_L - a L  \geq t +s )}{ \min_{K \leq N \leq L} 
    \PP (W_{L-N+1} - a (L-N)  \geq  s )} \, .
$$
\end {lemma}

\begin {proof}
Let $ B_K  = \{ W_K - a K  \geq t \}$ and, for $ K < N \leq L$,
$$
B_N  \, = \, \{ W_N - a N  \geq t \} \cap \bigcap_{K \leq M < N} \{ W_M - a M  < t \} .
$$
The sets $B_N$, $K \leq N \leq L$, are disjoint and
$$
\bigcup_{K \leq N \leq L} B_N \, = \, \Big \{ \max_{K \leq N \leq L} ( W_N - a N ) \geq t \Big \}.
$$
Then,
\beqs \begin {split}
\PP ( W_L - a L  \geq t +s )
   & \, \geq \, \sum_{K \leq N \leq L} \PP \big ( W_L - aL \geq t + s, B_N) \\
   & \, \geq \, \sum_{K \leq N \leq L} \PP \big ( W_{[N,L]} - a(L-N) \geq  s, B_N) \\
   & \, = \, \sum_{K \leq N \leq L} \PP \big ( W_{[N,L]} - a(L-N) \geq  s \big) \, \PP( B_N) \\
\end {split} \eeqs
where we successively used superadditivity and independence of $ W_{[N,L]}$ and $B_N$.
The conclusion follows since $W_{[N,L]}$ has the same distribution as $W_{L-N+1}$.
\end {proof}

We address the proof of the limsup theorem.
We argue similarly in the exponential and geometric cases, making clear which
tail inequality is used.

\begin {proof}[Proof of Theorem \ref {thm.limsup}]
Let $\phi : \N \to \R$ be defined by $\phi (n) = (\log \log n)^{2/3}$ if $ n \geq \eu^{\eu}$, 
and $\phi (n) = 1$ if not, and $ n _k = [\rho^k]$, $k \in \N$, for some $\rho >1$
to be made precise below.

We start with the upper bound. For $\beta >0$ and $k \geq 1$, let
$$
A_k \, = \, \bigg \{ \max_{n_{k-1} < N \leq n_k} \frac {{\widetilde H}_N}{\phi (N)} \geq \beta \bigg\} .
$$
We aim at showing that for every $\beta > \big (\frac 34 \big)^{2/3}$,
$\sum_k \PP (A_k) < \infty$, so that the conclusion follows by the Borel-Cantelli lemma.

By definition of ${\widetilde H}_N$, 
$$
\PP (A_k) \, \leq \, \PP \Big ( \max_{n_{k-1} < N \leq n_k} ( W_N - a N) \geq
            \beta b n_{k-1}^{1/3} \phi (n_{k-1}) \Big) .
$$
By the maximal inequality of Lemma~\ref {lem.maximal}, for any $s \geq 0$,
$$
\PP(A_k) \, \leq \, \frac {1}{D} \,
     \PP \big ( W_{n_k} - a n_k \geq \beta b n_{k-1}^{1/3} \phi (n_{k-1}) + s \big)
$$
where
$$
D \, = \,  \min_{n_{k-1} < N \leq n_k} \PP  \big ( W_{n_k - N +1} - a (n_k - N) \geq  s \big)
$$
For $ s = a$, \eqref {eq.rightlowerlr} ensures
that $ D \geq c >0$ independently of $k$. In the geometric case, we may rely
on \eqref {eq.rightbaik} for the choice of $ s = a + M b(n_k - N +1)^{1/3}$ for example.
Let then $\beta > \beta ' > \big (\frac 34)^{2/3}$.
Provided $\rho$ is close enough to $1$, for every $k$ large enough,
$$
\beta n_{k-1}^{1/3} \phi (n_{k-1}) \, \geq \, \beta' n_k^{1/3} \phi (n_k).
$$
Then, by \eqref {eq.rightupperj} (and \eqref {eq.j}), for every $\eta >0$ and
every $k$ large enough,
$$
\PP \big ( W_{n_k} \geq  a n_k + \beta b n_{k-1}^{1/3} \phi (n_{k-1})  \big)
   \, \leq \,  \eu^{- (\frac 43 -\eta)\beta'^{3/2}\phi (n_k)^{3/2}} .
$$
At this point therefore, for every $k$ large enough,
$$
\PP(A_k) \, \leq \, \frac {1}{c} \, \eu^{- (\frac 43 -\eta)\beta'^{3/2}\phi (n_k)^{3/2}} .
$$
Since $\beta ' > \big (\frac 34)^{2/3}$, there is $\eta >0$ such that the right-hand side of the
preceding inequality defines the general term of a convergent series.
Hence $\sum_k \PP (A_k) < \infty$ which completes the proof of the upper bound.

\medskip

Next, we turn to the lower bound. Recall that $ n _k = [\rho^k]$, $k \in \N$, 
where $\rho >1$. Assume first that there exists $\alpha >0$ such that for any $\rho >1$,
\beq \label {eq.infinite}
\sum_{k \geq 1}\PP \big ( W_{[n_{k-1},n_k]}  \geq a  (n_k - n_{k-1} + 1) 
     + \alpha b (n_k - n_{k-1} + 1)^{1/3} \phi (n_k - n_{k-1} + 1)\big) \, = \,  \infty.
\eeq 
By the independent part of the Borel-Cantelli lemma, on a set of probability one,
infinitely often in $k \geq1$,
$$
W_{[n_{k-1},n_k]}  \, \geq \,  a(n_k - n_{k-1} + 1) +
      \alpha b (n_k - n_{k-1} +1) ^{1/3} \phi (n_k - n_{k-1} + 1) .
$$
On the other hand, according to \eqref {eq.leftupperlr} in the exponential case
or \eqref {eq.leftbaik} in both the exponential and geometric cases, for any $\delta  >0$,
$$
\sum_{k \geq 1} \PP \big ( W_{n_{k-1}}  
     \leq  a n_{k-1} - \delta b  n_{k-1} ^{1/3} \phi ( n_{k-1}) \big)  \, < \, \infty.
$$
Hence, almost surely, for every $k$ large enough,
$$
W_{n_{k-1}}  \, \geq \, a n_{k-1} - \delta b  n_{k-1} ^{1/3} \phi ( n_{k-1}).
$$

As a consequence of the superadditivity inequality \eqref {eq.superadditivity}, 
on a set of probability one, infinitely often in $k$,
$$
W_{n_k}   \, \geq \,  an_k  +  \alpha b (n_k - n_{k-1}) ^{1/3} \phi (n_k - n_{k-1}) 
            - \delta b  n_{k-1} ^{1/3} \phi ( n_{k-1}).
$$
For every $\alpha ' < \alpha$, if $\rho >1$ is large enough,
$$
\alpha b (n_k - n_{k-1}) ^{1/3} \phi (n_k - n_{k-1}) 
            - \delta b  n_{k-1} ^{1/3} \phi ( n_{k-1})
            \, \geq \,  \alpha' b n_k  ^{1/3} \phi (n_k ) .
$$
Hence, since $\alpha' < \alpha$ is arbitrary,
$$
\liminf_{N \to \infty} \frac {{\widetilde H}_N}{\phi (N)} \, = \,
\liminf_{N \to \infty} \frac {W_N - a N}{bN^{1/3} \phi (N) } \, \geq \, \alpha
$$
almost surely.

It remains to discuss the choice of $\alpha >0$ so that \eqref {eq.infinite} holds.
Set $m_k = n_k - n_{k-1} + 1$.
On the basis of \eqref {eq.rightlowerlr}, for some $c ,C>0$ and every $k \geq 1$ large enough,
$$
\PP \big ( W_{m_k}  \geq  a  m_k + \alpha b m_k^{1/3} \phi (m_k)\big)
   \, \geq \, c \,  \eu^{-  C \alpha ^{3/2} \phi (m_k)^{3/2}  } .
$$
Provided $\alpha >0 $ is small enough, \eqref {eq.infinite} is satisfied.
Now, if we agree that \eqref {eq.rightbaik} holds true, for some $C >0$ and every $k \geq 1$ large enough,
$$
\PP \big ( W_{m_k}  \geq  a  m_k + \alpha b m_k^{1/3} \phi (m_k)\big)
   \, \geq \, \phi(m_k)^{-C} \,  \eu^{-  \frac 43 \alpha ^{3/2} \phi (m_k)^{3/2}  } .
$$
In this case, \eqref {eq.infinite} is satisfied for all $\alpha < \big (\frac 34 \big)^{2/3}$,
yielding the conjectured lower bound in Theorem~\ref {thm.limsup}.
\end {proof}

Next, we turn to the liminf theorem. Since the superadditivity property is only one-sided,
a different (weaker) strategy has to be followed,
yielding in particular non-optimal bounds.

\begin {proof}[Proof of Theorem \ref {thm.liminf}]
Let $\psi (n) = (\log \log n)^{1/3}$ if $ n \geq \eu^{\eu}$, and $\psi (n) = 1$ if not.
Let $ 0 < \eta < 1$ and set here $n_k = [\eu^{k^\eta}]$, $k \geq 1$. 

By the Borel-Cantelli lemma, it is enough to establish that $\sum_k \PP (A_k) < \infty$ where
$$
A_k \, = \, \bigg \{ \min_{n_{k-1} < N \leq n_k} \frac {{\widetilde H}_N}{\psi (N)} \leq - 2\alpha \bigg\}
$$
for some (large enough) $ \alpha >0$. 
For every $k \geq 1$,
$$
\PP (A_k) 
      \, \leq \, \sum_{N = n_{k-1}+1}^{n_k} 
        \PP \bigg ( \frac {{\widetilde H}_N}{\psi (N)} \leq - 2\alpha, 
           \frac {{\widetilde H}_{n_{k-1}}}{\psi (n_{k-1})} \geq -  \alpha \bigg)
        +  \PP \bigg ( \frac {{\widetilde H}_{n_{k-1}}}{\psi (n_{k-1})} \leq -  \alpha \bigg).
$$
Now
$$
\PP \bigg ( \frac {{\widetilde H}_{n_{k-1}}}{\psi (n_{k-1})}  \leq -  \alpha \bigg)
   \, = \, \PP \big ( W_{n_{k-1}} \leq (a - \varepsilon ) n_{k-1}\Big)
$$
where $\varepsilon n_{k-1} = \alpha b n_{k-1}^{1/3} \,\psi (n_{k-1})$.
By \eqref {eq.leftupperlr} in the exponential case
or \eqref {eq.leftbaik} in both the exponential and geometric cases, 
$$
\PP \big ( W_{n_{k-1}} \leq (a - \varepsilon ) n_{k-1}\Big)
   \, \leq \, C \, \eu^{- c (\alpha b)^3 \psi (n_{k-1})^3} .
$$
The right-hand side defines the general term of a convergent series whenever $\alpha >0$ is large enough.

Next, by superadditivity \eqref {eq.superadditivity},
$$
\PP \bigg ( \frac {{\widetilde H}_N}{\phi (N)}  \leq - 2\alpha ,
               \frac {{\widetilde H}_{n_{k-1}}}{\phi (n_{k-1})} \geq - \alpha \bigg) \\
      \, \leq \, \PP \big ( W_{ N - n_{k-1}+1}  \leq    (a - \varepsilon )(  N - n_{k-1}) \big) 
$$
where now
$$
\varepsilon (N - n_{k-1} )
     \, = \, \alpha b \big [2  N^{1/3} \psi (N) -  n_{k-1}^{1/3} \,\psi (n_{k-1}) \big ]
$$
(assumed to be strictly positive).  By \eqref {eq.leftupperlr} or \eqref {eq.leftbaik} again,
$$
\PP \big ( W_{ N - n_{k-1} + 1}  \leq    (a - \varepsilon )(  N - n_{k-1}) \big) 
            \, \leq \, C \eu^{-c \varepsilon ^3 ( N - n_{k-1})^2} .
$$
Now, for every $n_{k-1} < N \leq n_k$,
$$
2  N^{1/3} \psi (N) -  n_{k-1}^{1/3} \,\psi (n_{k-1}) \, \geq \,  N^{1/3} 
$$
so that $\varepsilon (N - n_{k-1} ) \geq \alpha b N^{1/3}$.
In addition, for some $\delta >0$ and every $k$ large enough,
$$
\frac {N}{N- n_{k-1}} \, \geq \, \delta \, k^{1-\eta} \, .
$$
Hence,
$$
 \sum_{N = n_{k-1}+1}^{n_k} 
   \eu^{-c \varepsilon ^3 (N - n_{k-1})^2}
   \, \leq \, \sum_{N = n_{k-1}+1}^{n_k}  \eu^{- c \delta (\alpha b)^3  k^{1 - \eta}}
   \, \leq \,  \eu^{k^\eta} \eu^{- c \delta (\alpha b)^3  k^{1-\eta}} .
$$
Provided $\eta >0$ is small enough, the right-hand side defines the general
term of convergent series in $k$. Together with the previous step,
$\sum_k \PP(A_k) < \infty$, and the proof of Theorem~\ref {thm.liminf} is complete.
\end {proof}

\vskip 8mm

\font\tenrm =cmr10  {\tenrm

\parskip 0mm

\noindent Institut de Math\'ematiques de Toulouse 

\noindent Universit\'e de Toulouse -- Paul-Sabatier, F-31062 Toulouse, France

\noindent \&  Institut Universitaire de France 

\noindent ledoux@math.univ-toulouse.fr

}

\end {document}